\documentclass[11pt,a4paper]{article}
\usepackage{amssymb}
\usepackage{amsmath, amsthm, amscd, amsfonts, amssymb, graphicx, color}
\usepackage[bookmarksnumbered, colorlinks, plainpages]{hyperref}
\newtheorem{thm}{Theorem}[section]

\newtheorem{lem}[thm]{Lemma}

\newtheorem{obs}[thm]{Observation}

\newtheorem{prop}[thm]{Proposition}
\newtheorem{defn}[thm]{Definition}

\numberwithin{equation}{section}
\numberwithin{equation}{subsection}

\begin{document}

\title{Total Dominator Total Chromatic Numbers of Wheels, Complete bipartite graphs and Complete graphs}

\author{$^{1}$Adel P. Kazemi and $^{2}$Farshad Kazemnejad% and $^{3}$Somayeh Moradi
 \\[1em]
$^{1,2}$ Department of Mathematics, \\ University of Mohaghegh Ardabili, \\ P.O.\ Box 5619911367, Ardabil, Iran. \\
$^1$ Email: adelpkazemi@yahoo.com\\
$^2$ Email:  kazemnejad.farshad@gmail.com\\
[1em]
}

\maketitle

\begin{abstract}
Total dominator total coloring of a graph is a total coloring of the graph such that each object of the graph is adjacent or incident to every object of some color class. The minimum namber of the color classes of a total dominator total coloring of a graph is called the total dominator total chromatic number of the graph. Here, we will find  the total dominator chromatic numbers of  wheels, complete bipartite graphs and complete graphs.
\\[0.2em]

\noindent
Keywords: Total dominator total coloring, Total dominator total chromatic number, total domination number, total mixed domination number, total graph.
\\[0.2em]

\noindent
MSC(2010): 05C15, 05C69.
\end{abstract}

%------------------------------------------------------------------------------------%
\pagestyle{myheadings}
\markboth{\centerline {\scriptsize A. P. Kazemi and F. Kazemnejad}}     {\centerline {\scriptsize A. P. Kazemi and F. Kazemnejad,~~~~~~~~~~~~~~~TDTCN of wheels, complete bipartite graphs and complete graphs}}

%------------------------------------------------------------------------------------

\section{\bf Introduction}

All graphs considered here are non-empty, finite, undirected and simple. For
standard graph theory terminology not given here we refer to \cite{West}. Let $%
G=(V,E) $ be a graph with the \emph{vertex set} $V$ of \emph{order}
$n$ and the \emph{edge set} $E$ of \emph{size} $m$. The
\emph{open neighborhood} and the \emph{closed neighborhood} of a
vertex $v\in V$ are $N(v)=N_{G}(v)=\{u\in V\ |\ uv\in E\}$ and
$N[v]=N_{G}[v]=N_{G}(v)\cup \{v\}$, respectively. The \emph{degree} of a
vertex $v$ is also $deg_G(v)=\vert N_{G}(v) \vert $. The
\emph{minimum} and \emph{maximum degree} of $G$ are denoted by
$\delta =\delta (G)$ and $\Delta =\Delta (G)$, respectively. If
$\delta (G)=\Delta (G)=k$, then $G$ is called $k$-\emph{regular}. An \emph{independent set} of $G$ is a subset of vertices of $G$, no two of which are adjacent. And a \emph{maximum independent set} is an independent set of the largest cardinality in $G$. This cardinality is called the \emph{independence number} of $G$, and is denoted by $\alpha(G)$. Also a \emph{mixed independent set} of $G$ is a subset of $V\cup E$, no two objects of which are adjacent or incident, and a \emph{maximum mixed  independent set} is a mixed independent set of the largest cardinality in $G$. This cardinality is called the \emph{mixed independence number} of $G$, and is denoted by $\alpha_{mix}(G)$. Two isomorphic graphs $G$ and $H$ are shown by $G\cong H$. We write $K_{n}$ , $C_{n}$ 
and $P_{n}$ for a \emph{complete graph}, a \emph{cycle} and a \emph{path} 
of order $n$, respectively, while $W_{n}$, $K_{m,n}$ and  $G[S]$ denote a \emph{wheel} of order $n+1$,   a \emph{complete bipartite graph} of order $m+n$ and   the \emph{induced subgraph} of $G$ by a vertex set $S$, respectively.

\vskip 0.2 true cm

\begin{figure}[ht]
\centerline{\includegraphics[width=9cm, height=2.2cm]{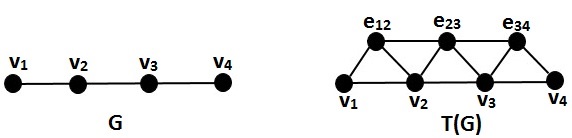}}
\vspace*{-0.3cm}
\caption{ The illustration of $G$ (left) and $T(G)$ (right).}\label{TDTC100}
\end{figure}

 The \emph{Cartesian product} $G \square H$ of two graphs $G$ and $H$ is a graph with $V(G) \times V(H)$ and two vertices  $(g_{1}, h_{1})$ and $(g_{2}, h_{2})$ are adjacent if and only if either $g_{1}=g_{2}$ and $(h_{1}, h_{2}) \in E(H)$, or $h_{1}=h_{2}$ and $(g_{1}, g_{2}) \in E(G)$. The \emph{line graph} $L(G)$ of $G$ is a graph with the vertex set $E$ and two vertices of $L(G)$ are adjacent when they are incident in $G$.  The \emph{total graph} $T(G)$ of a graph $G$ is the graph whose vertex set is $V\cup E$ and two vertices are adjacent whenever they are either adjacent or incident in $G$ \cite{Behzad}. It is obvious that if $G$ has order $n$ and size $m$, then $T(G)$ has order $n+m$ and size $3m+|E(L(G))|$, and also $T(G)$ contains both $G$ and $L(G)$ as two induced subgraphs and it is the largest graph formed by adjacent and incidence relation between graph elements. In this paper, by assumption $V=\{v_1,v_2,\cdots, v_n\}$, we use the notations $V(T(G))=V\cup \mathcal{E}$ where $\mathcal{E}=\{e_{ij}~|~v_iv_j\in E\}$, and $E(T(G))=\{v_ie_{ij},v_je_{ij}~|~v_iv_j\in E\}\cup E\cup E(L(G))$. Obviousely $deg_{T(G)}(v_i)=2deg_G(v_i)$ and $deg_{T(G)}(e_{ij})=deg_G(v_i)+deg_G(v_j)$. So if $G$ is $k$-regular, then $T(G)$ is $2k$-regular. Also $\alpha_{mix}(G)=\alpha(T(G))$. In Figure \ref{TDTC100} a graph $G$ and its total graph are shown for an example.

%-------------------------------------------------------------
\vskip 0.2 true cm

\textbf{DOMINATION.} Domination in graphs is now well studied in graph theory and the literature
on this subject has been surveyed and detailed in the two books by Haynes,
Hedetniemi, and Slater~\cite{hhs1, hhs2}. A famous type of domination is total domination, and the literature on this subject has been surveyed and detailed in the recent
book~\cite{HeYe13}. A \emph{total dominating set}, briefly TDS, $S$ of a graph $G=(V,E)$ is a subset
of the vertex set of $G$ such that for each vertex $v$, $N_G(v)\cap S\neq \emptyset$. The \emph{total domination number $\gamma_{t}(G)$} of $G$ is the minimum cardinality of a TDS of $G$. Similarly, a subset $S\subseteq V\cup E$ of a graph $G$ is called a \emph{total mixed dominating set}, briefly TMDS, of $G$ if each object of $V\cup E$ is either adjacent or incident to an object of $S$, and the \emph{total mixed domination number} $\gamma_{tm}(G)$ of $G$ is the minimum cardinality of a TMDS \cite{KK2017}. A min-TDS/min-TMDS of $G$ denotes a TDS/TMDS of $G$ with minimum cardinality. Also we agree that \emph{a vertex $v$ dominates an edge} $e$ or \emph{an edge $e$ dominates a vertex} $v$ mean $v\in e$. Similarly, we agree that \emph{an edge dominates another edge} means they have a common vertex.  The next theorem can be easily obtained.

\begin{thm} \emph{\cite{KK2017}}
\label{gamma_{tm}(G)=gamma_{t}(T(G))}
For any graph $G$ without isolate vertex, $\gamma_{tm}(G)=\gamma_{t}(T(G))$.
\end{thm}
%----------------------------------------------------------------
\vskip 0.2 true cm
\textbf{GRAPH COLORING.} Graph coloring is used as a model for a vast
number of practical problems involving allocation of scarce
resources (e.g., scheduling problems), and has played a key role in
the development of graph theory and, more generally, discrete
mathematics and combinatorial optimization. Graph colorability
is NP-complete in the general case, although the problem is solvable
in polynomial time for many classes \cite{GJ}. A \emph{proper coloring} of a graph $G$ is a function from
the vertices of the graph to a set of colors such that any two adjacent vertices have different colors, and the minimum number of colors needed in a proper
coloring of a graph is called the \emph{chromatic number} $\chi (G)$ of $G$. In a simlar way, a \emph{total coloring} of $G$ assigns a color to each vertex and to each edge so that colored objects have different colors when they are adjacent or incident, and the minimum number of colors needed in a total coloring of a graph is called the \emph{total chromatic number} $\chi_{T}(G)$ of $G$ \cite{West}. %The Total Coloring Conjecture (Behzad [1965]) states that 
%\begin{BehConj}
%For every simple graph $G$, $\chi_{T}(G) \leq \Delta(G) + 2$. 
%\end{BehConj}
 A \emph{color class} in a coloring of a graph is a set consisting of all those objects assigned the same color. For simply, if $f$ is a (total) coloring of $G$ with the coloring classes $V_1$, $V_2$, $\cdots$ , $V_{\ell}$, we write $f=(V_1,V_2,\cdots,V_{\ell})$. Hence $V=V_1 \cup V_2\cup \cdots \cup V_{\ell}$ is a partition of the vertex set of the graph, and so
%------------------------
 \begin{equation}
 \label{|V|=sum_{i=1}^{ell} |V_i|}
 |V|=\sum_{i=1}^{\ell} |V_i|.
\end{equation}
 %--------------------
 Motivated by the relation between coloring and total dominating, the concept of total dominator coloring in graphs introduced in \cite{Kaz2015} by Kazemi, and extended in \cite{Hen2015,Kaz-Par,Kaz2014,Kaz2016,KK2018}.

%coloring of the edges of a nonempty graph is defined by Chartrand and et al. in \cite{GP} and $\chi '(G)$ denotes the \emph{edge chromatic number} (or \emph{chromatic index}) of $G$. Trivially $\chi^{'}(G)=\chi(L(G))$.
%-----------------------------------------------
%\vskip 0.2 true cm
%---------------------------------------------------------------------------------------------------------------
\begin{defn} 
\label{total dominator coloring} \emph{ \cite{Kaz2015} A} total dominator coloring,
\emph{briefly TDC, of a graph $G$ with a possitive minimum degree is a proper coloring of $G$ in
which each vertex of the graph is adjacent to every vertex of some
color class. The }total dominator chromatic number $\chi_{d}^t(G)$
\emph{of $G$ is the minimum number of color classes in a TDC of $G$.}
\end{defn}

%----------------------------------------------------------
In \cite{KKM2019}, the authors initiated studying of a new concept called total dominator total coloring in graphs which is obtained from the concept of total dominator coloring of a graph by replacing total coloring of a graph instead of (vertex) coloring of it.
%----------------------------------------------------------
\begin{defn} 
\label{total dominator total coloring} \emph{ A} total dominator total coloring,
\emph{briefly TDTC, of a graph $G$ with a possitive minimum degree is a total coloring of $G$ in
which each object of the graph is adjacent or incident to every object of some
color class. The }total dominator total chromatic number $\chi_d^{tt}(G)$
\emph{of $G$ is the minimum number of color classes in a TDTC of $G$.}
\end{defn}
%---------------------------------------------------------------------------------------------------------------
%\textcolor{blue}{\begin{defn}\label{free color class} \emph{\cite{Hen2015}}\emph{A color class in a TDC $\mathcal{C}$ of $G$ is called} free \emph{if each vertex of $G$ is adjacent to every vertex of some color class different from it.(((By Farshad)))} \end{defn}}

It can be easily obtained the next theorem.

\begin{figure}[ht]\label{TDTC11}
\centerline{\includegraphics[width=10cm, height=2.7cm]{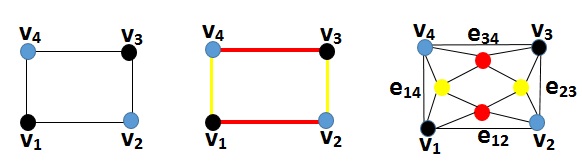}}
\vspace*{-0.1cm}
\caption{A min-TDC of $C_4$ (left), a min-TDTC of $C_4$ (Middle) and its corresponding min-TDC of $T(C_4)$ (right).} \label{TDTC101}
\end{figure}

\begin{thm}
\emph{\cite{KKM2019}}
\label{chi_d^{tt}(G)=chi_d^{t}(T(G))}
For any graph $G$ without isolate vertex, $\chi_d^{tt}(G)=\chi_d^{t}(T(G))$.
\end{thm}
%----------------------------------------------------------
For any TDC (TDTC) $f=(V_1,V_2,\cdots,V_{\ell})$ of a graph $G$, a vertex (an object) $v$ is called a \emph{common neighbor} of $V_i$ or we say $V_i$ \emph{totally dominates} $v$, and we write $v\succ_t V_i$, if vertex (object) $v$ is adjacent (adjacent or incident) to every vertex (object) in $V_i$. Otherwise we write $v \not\succ_t  V_i$. Also $v$ is called a \emph{private neighbor of $V_i$ with respect to $f$} if $v\succ_t V_i$ and $v\nsucc_t V_j$ for all $j\neq i$. The set of all common neighbors of $V_i$ with respect to $f$ is called the \emph{common neighborhood} of $V_i$ in $G$ and denoted by $CN_{G,f}(V_i)$ or simply by $CN(V_i)$. Also every TDC or TDTC of $G$  with $\chi_d^t(G)$ or $\chi_d^{tt}(G)$ colors is called respectively a \emph{min}-TDC or a \emph{min}-TDTC. For an example see Figure \ref{TDTC11}.
%---------------------------------------------------
%------------------------------------------------------

Also for any TDC $(V_1,V_2,\cdots , V_{\ell})$ and any TDTC $(W_1,W_2,\cdots , W_{\ell})$  of a graph $G$, we have
\begin{equation}
\label{CN(V_i)=V}
\bigcup_{i=1}^{\ell} CN(V_i)=V \mbox{ and } \bigcup_{i=1}^{\ell} CN(W_i)=V\cup E.
\end{equation}

\vspace{0.2cm}

\textbf{GOAL.} In \cite{KKM2019}, the authors initiated to study the total dominator total coloring of a graph and found some useful results, and presented some problems. Finding the total dominator total chromatic numbers of wheels, complete bipartite graphs and complete graphs were three of them which we consider them here. We recall the following propositions which are useful for our investigation. 
%-------------------------------------------------------------
\begin{prop}\emph{\cite{Kaz2015}}
\label{chi_d^t =<g_t+min{chi(G-S)}} For any connected graph $G$ with $\delta(G)\geq 1$, 
%-------------------------
\begin{equation*}
\label{chi_d^t(G) =<gama_t(G)+ min chi(G[V(G)-S])}
\chi_d^t(G) \leq \gamma_{t}(G)+ \min_S \chi(G[V(G)-S]),
\end{equation*}
%---------------------------
where $S\subseteq V(G)$ is a min-TDS of $G$. And so $\chi_d^t(G) \leq \gamma_t(G)+ \chi(G)$.
\end{prop}
%------------------------------------------------------------------------
\begin{prop}\emph{\cite{KK2019}}
\label{alpha (T(P_n))=lfloor 2n-1/3 rfloor}
For any path $P_n$ of order $n\geq 3$, $\alpha_{mix}(P_n)=\lceil \frac{2n-1}{3}\rceil$.
\end{prop}
%--------------------------------------------------------------
\begin{prop}\emph{\cite{KK2019}}
\label{alpha (T(C_n))=lfloor 2n/3 rfloor}
For any cycle $C_n$ of order $n\geq 3$, $\alpha_{mix}(C_n)=\lfloor \frac{2n}{3}\rfloor$.
\end{prop}
%--------------------------------------------------------
\begin{prop}\emph{\cite{KK2019}}
\label{chi^{t}_d T(P_n)}
For any path $P_n$ of order $n\geq 2$,
\begin{equation*}
\chi^{tt}_d (P_n)=\left\{
\begin{array}{ll}
\gamma_{tm}(P_n)+1  & \mbox{if }n=2,3, \\
\gamma_{tm}(P_n)+2  & \mbox{if }n=4,5,6,8,9,10,13,16, \\
\gamma_{tm}(P_n)+3  & \mbox{if }n= 7,\mbox{ } n\neq 13,16\mbox{ or }n \geq 11 
\end{array}
\right.
\end{equation*}
%----------------
or equivalently
%---------------
\begin{equation*}
\chi^{tt}_d (P_n)=\left\{
\begin{array}{ll}
n+1 & \mbox{if }n=2,\\
n     & \mbox{if }3\leq n \leq 7,\\
n-1  & \mbox{if }8\leq n \leq 9,%\\
%n-2=8     & \mbox{if }n=10,
\end{array}
\right.
\end{equation*}
and for $n \geq 10$,
\begin{equation*}
\chi^{tt}_d (P_n)=\left\{
\begin{array}{ll}
\lfloor \frac{4n}{7}\rfloor +3 & \mbox{if }n\equiv 4 \pmod{7} \mbox{ or } n=10, 13, 16,\\
\lceil \frac{4n}{7}\rceil   +3   & \mbox{if }n\not\equiv 4 \pmod{7} \mbox{ and } n\neq 10, 13,16.
\end{array}\right.
\end{equation*}
\end{prop}
%---------------------------------------------------
\begin{prop}\emph{\cite{KK2019}}
\label{chi^{t}_d T(C_n)}
For any cycle $C_n$ of order $n\geq 3$,
\begin{equation*}
\chi^{tt}_d(C_n)=\left\{
\begin{array}{ll}
\gamma_{tm}(C_n)+1  & \mbox{if }n=3,4,5 \\
\gamma_{tm}(C_n)+2  & \mbox{if }n=6,9,12 \\
\gamma_{tm}(C_n)+3  & \mbox{if }n\geq 7 ~\mbox{and}~n \neq 9,12.
\end{array}
\right.
\end{equation*}
%-----------------
or equivalently
%-------
\begin{equation*}
\chi^{tt}_d(C_n)=\left\{
\begin{array}{ll}
n & \mbox{if }3\leq n \leq 8,\\
n-1 & \mbox{if }n=9,
\end{array}\right.
\end{equation*}
%---------------
and for $n\geq 10$,
%-----------------------
\begin{equation*}
\chi^{tt}_d(C_n)=\left\{
\begin{array}{ll}
\lceil \frac{4n}{7}\rceil +4 & \mbox{if }n\equiv 5 \pmod{7} \mbox{ and } n\neq 12,\\
\lceil \frac{4n}{7}\rceil +3 & \mbox{if }n\not\equiv 5 \pmod{7} \mbox{ or } n=12.
\end{array}\right.
\end{equation*}
%----------------
\end{prop}
%---------------------------------------------------------
\section{Wheels}

Here, we calculate the total dominator total chromatic number of a wheel. First we  recall a proposition from \cite{KK2017} and calculate the mixed indepence number of a wheel.

\begin{prop} \emph{\cite{KK2017}}
\label{gamma_t(T(W_n))=lceil n/2 rceil +1}  For any wheel $W_n$ of order $n+1\geq 4$, $\gamma_{tm}(W_n)=\lceil \frac{n}{2}  \rceil +1$.
\end{prop}

%------------------------------------------------------------
\begin{lem}\label{alpha (T(W_n))=lceil 2n/3 rceil}
For any wheel $W_n$ of order $n+1\geq 4$, $\alpha_{mix}(W_n)=\lceil \frac{2n}{3}\rceil$.
\end{lem}

\begin{proof} 
Let  $W_n=(V,E)$ be a wheel of order $n+1\geq 4$ where $ V=\{v_{i}~|~ 0 \leq i \leq n \} $ and $ E=\{v_{0}v_{i}, v_{i}v_{i+1}~|~ 1\leq i \leq n \}$. 
Then $ V(T(W_{n}))=V\cup \mathcal{E}$ when $\mathcal{E}=\{e_{0i}, e_{i(i+1)}~|~1\leq i \leq n\}$.  Let $S$ be an independent set of $T(W_n)$. Since the subgraph induced by $\{e_{0i}~|~1\leq i \leq n\}\cup\{v_0\}$ is a complete graph,  we have $|S\cap (\{e_{0i}~|~1\leq i \leq n\}\cup\{v_0\})|\leq 1$.   If $|S\cap (\{e_{0i}~|~1\leq i \leq n\}\cup\{v_0\})|=0$, then $S\subseteq \{v_i,e_{i(i+1)}~|~1 \leq i \leq n\}$, and since the subgraph induced by $\{v_i,e_{i(i+1)}~|~1 \leq i \leq n\}$ is isomorphic to $T(C_n)$, Proposition \ref{alpha (T(C_n))=lfloor 2n/3 rfloor} implies $|S|\leq \lfloor \frac{2n}{3}\rfloor$. If also $v_0  \in S$, then $S\subseteq \{e_{i(i+1)}~|~1 \leq i \leq n\}$, and since the subgraph induced by $\{e_{i(i+1)}~|~1 \leq i \leq n\}$ is isomorphic to $C_n$, we have $|S|\leq \alpha(C_n)+1= \lfloor \frac{n}{2}\rfloor+ 1$. Finally if $e_{0i} \in S$ for some $1 \leq i \leq n$, then $S\subseteq V\cup \mathcal{E}-N_{T(W_n)}(e_{0i})$, and since the subgraph induced by $V\cup \mathcal{E}-N_{T(W_n)}(e_{0i})$ is isomorphic to $T(P_{n-1})$, Proposition \ref{alpha (T(P_n))=lfloor 2n-1/3 rfloor} implies $|S|\leq \lceil \frac{2n}{3}\rceil$. Therefore $\alpha_{mix}(W_n)=\alpha (T(W_n))=\max\{\lfloor \frac{2n}{3}\rfloor, \lfloor \frac{n}{2}\rfloor+ 1,\lceil \frac{2n}{3}\rceil \}=\lceil \frac{2n}{3}\rceil$. 
\end{proof}

%-------------------------------------------------------------------------------------------
\begin{prop}
\label{chi_d^t(T({W_n}))}
 For any wheel $W_n$ of order $n+1\geq 4$,
 \begin{equation*}
\chi_d^{tt}(W_n)=\left\{
\begin{array}{ll}
n+2   & \mbox{if }3 \leq n \leq 7,\\
n+1 & \mbox{if }n\geq 8.
\end{array}
\right.
\end{equation*}
\end{prop}
\begin{proof}
Let  $W_n=(V,E)$ be a wheel of order $n+1\geq 4$ where $ V=\{v_{i}~|~ 0 \leq i \leq n \} $ and $ E=\{v_{0}v_{i}, v_{i}v_{i+1}~|~ 1\leq i \leq n \}$. 
Then $ V(T(W_{n}))=V\cup \mathcal{E}_0 \cup \mathcal{E}_1$ when $\mathcal{E}_0=\{ e_{0i}~|~1\leq i \leq n\}$ and $\mathcal{E}_1=\{e_{i(i+1)}~|~1\leq i \leq n\}$.  Let $ f=(V_{1},\cdots,V_{\ell})$ be a min-TDC of $T(W_{n})$. Since the subgraph of $T(W_n)$ induced by $\mathcal{E}_0\cup \{v_0\}$ is isomorphic to a complete graph of order $n+1$, we have
\begin{equation}\label{chi_d^t(T({W_n})) geq n+1}
\chi_d^t(T({W_n}))\geq n+1.
\end{equation}
%-------------------------
For $n\geq8$, by  Proposition \ref{chi_d^t =<g_t+min{chi(G-S)}}, we know
\[
\chi_d^{tt}(W_n) \leq \gamma_t(T(W_n))+\chi(G[V(T(W_n))-S]),
\]
when $S$ is a min-TDS of $T({W_n})$. Since, by  Proposition \ref{gamma_t(T(W_n))=lceil n/2 rceil +1}, the sets $S_e=\{e_{0(2i)}~|~1\leq i \leq \lfloor \frac{n}{2} \rfloor\}\cup \{v_{0}\}$, when $n$ is even, and $S_o=\{e_{0(2i)}~|~1\leq i \leq \lfloor \frac{n}{2}  \rfloor\}\cup \{v_{0}, e_{0n}\}$, when $n$ is odd, are two min-TDSs of $T({W_n})$  of cardinality $\lceil \frac{n}{2} \rceil+1$, we have
\[
\chi_d^{tt}(W_n)  \leq \lceil \frac{n}{2} \rceil+1+\chi(G[V(T(W_n))-S]),
\]
in which $S=S_e$ when $n$ is even and $S=S_o$ when $n$ is odd. To complete our proof it is sufficient to prove $\chi(G[V(T(W_n))-S])= \lfloor \frac{n}{2} \rfloor$. Since the  subgraph induced by $\{e_{0(2i-1)}~|~\ 1\leq i\leq \lfloor \frac{n}{2} \rfloor\}$ is a complete graph, we have $\chi(G[V(T(W_n))-S]) \geq \lfloor \frac{n}{2} \rfloor$. On the other hand, since, for even $n$ the coloring function $f_e$ with the criterion
 \begin{equation*}
f_e(w) \equiv \left\{
\begin{array}{ll}
i ~~~~~  \pmod {\lfloor {\frac{n}{2}} \rfloor}       & \mbox{if }w=e_{0(2i+1)},\\
i+1 \pmod {\lfloor {\frac{n}{2}} \rfloor}   & \mbox{if }w=e_{(2i+1)(2i+2)} \mbox{ or }v_{2i+3},\\
i+2 \pmod {\lfloor {\frac{n}{2}} \rfloor}   & \mbox{if }w=v_{2i+2},\\
i+3 \pmod {\lfloor {\frac{n}{2}} \rfloor}   & \mbox{if }w=e_{(2i+2)(2i+3)}
\end{array}
\right.
\end{equation*}
when $0 \leq i \leq \lfloor {\frac{n}{2}} \rfloor-1$ is a proper coloring of $G[V(T(W_n))-S_e$, and for odd $n$ the coloring function $f_o$ with the criterion
 \begin{equation*}
f_o(w) \equiv \left\{
\begin{array}{ll}
i ~~~~~\pmod {\lfloor {\frac{n}{2}} \rfloor}                 & \mbox{if }w=e_{0(2i+1)},\\
i+1 \pmod {\lfloor {\frac{n}{2}} \rfloor}                         & \mbox{if }w=e_{(2i+1)(2i+2)} \mbox{ or }v_{2i+3},\\
i+2 \pmod {\lfloor {\frac{n}{2}} \rfloor}                         & \mbox{if }w=v_{2i+2},\\
i+3 \pmod {\lfloor {\frac{n}{2}}\rfloor}                          & \mbox{if }w=e_{(2i+2)(2i+3)},\\
2                                                                                  & \mbox{if }w=e_{(n-1)n},\\
3                                                                                  & \mbox{if }w=v_n
\end{array}
\right.
\end{equation*}
when $0 \leq i \leq \lfloor {\frac{n}{2}} \rfloor-1 $, is a proper coloring of $G[V(T(W_n))-S_o$, we have $\chi(G[V(T(W_n))-S])= \lfloor \frac{n}{2} \rfloor$. Therefore, we continue our proof when $3 \leq n \leq 7$ by considering the following facts in which $f=(V_{1},V_{2},\cdots,V_{\ell})$ is a min-TDC of $T(W_{n})$, $|V_1|\geq |V_2|\geq \cdots \geq |V_{\ell}|$, and $\mathcal{B}_{i}=\{V_{k}~|~e_{i(i+1)} \succ_{t} V_{k}\mbox{ and }|V_{k}|= i \mbox{ for some } e_{i(i+1)}\in \mathcal{E}_1\}$, $b_i=|\mathcal{B}_i|$ for $1 \leq i \leq \lceil \frac{2n}{3} \rceil $.
% Also, without loss of generality, we assume $t_{0i} \in V_{i}$ for $1\leq i \leq n$, $v_{0} \in V_{n+1}$.
 %-------------------------
\begin{itemize}

\item[$\star$] \texttt{Fact 1.} $\sum_{i=1}^{\ell}|V_i|=3n+1$, by (\ref{|V|=sum_{i=1}^{ell} |V_i|}), and $3n+1\leq \ell \lceil \frac{2n}{3}\rceil$.% by $|V_{k}|\leq \alpha$ for each $1\leq k\leq \ell$.
 
\item[$\star$] \texttt{Fact 2.} For any $v \in \mathcal{E}_0\cup \mathcal{E}_1$, if $v \succ_{t} V_{k}$ for some $1\leq k\leq \ell$, then $|V_{k}|\leq 2$.

\item[$\star$] \texttt{Fact 3.} If $e_{i(i+1)} \succ_{t} V_{k}$ for some $e_{i(i+1)}\in \mathcal{E}_1$ and some $1\leq k\leq \ell$ and $|V_k|=2$, then $CN(V_k)\cap \mathcal{E}_1=\{e_{i(i+1)}\}$.

\item[$\star$] \texttt{Fact 4.} If $e_{i(i+1)} \succ_{t} V_{k}$ for some $1\leq k\leq \ell$ and $|V_k|=1$, then $|CN(V_k)\cap \mathcal{E}_1|=2$. %$V_k$ totally dominates only other one vertex  in $\mathcal{E}_1$.
 
\item[$\star$] \texttt{Fact 5.} $n \leq 2b_1+b_2\leq \ell$ (by \texttt{Facts 3, 4}).
 
\item[$\star$] \texttt{Fact 6.} For $1\leq i \leq n$, if $v_i \succ_{t} V_{k}$ for some $1\leq k\leq \ell$, then $|V_{k}|\leq 3$.

\item[$\star$] \texttt{Fact 7.} For $1\leq i \leq n$, if $v_i \succ_{t} V_{k}$ for some $1\leq k\leq \ell$ and $|V_{k}|=3$, then $CN(V_k)=\{v_0\}$.
 
\item[$\star$] \texttt{Fact 8.}  If $v_0 \succ_{t} V_{k}$ for some $1\leq k\leq \ell$, then $|V_{k}|\leq \lfloor \frac{n}{2} \rfloor +1$.
\end{itemize}

\begin{itemize}

\item $n=3$. Then $|V_i|\leq \alpha= 2$ for each $i$, and so $\ell \geq 5$, by \texttt{Fact 1}. Now since the coloring function $(\{e_{12},e_{03}\},\{v_{1},e_{23}\},\{v_{0},e_{13}\},\{v_{2},e_{01}\},\{v_{3},e_{02}\})$ is a TDC of $T(W_{3})$, we have $\chi_d^{tt}({W_3})=5$.
%------------------------------------------------------------
\item $n=4$. Then $|V_i|\leq \alpha= 3$ for each $i$, and so $\ell \geq 5$, by \texttt{Fact 1}. If $\ell=5$,  then $(|V_1|,|V_2|, \cdots ,|V_{5}|)=(3,3,3,3,1)$ which contradicts the \texttt{Facts 2, 4}, or $(|V_1|,|V_2|, \cdots ,|V_{5}|)=(3,3,3,2,2)$ which contradicts the \texttt{Facts 2, 3}. So $\ell \geq 6$. Now since $(V_1,\cdots,V_{6})$ is a TDC of $T(W_{4})$ where $V_{i}=\{e_{0i},v_{i+1}\}$ for $1\leq i \leq 3$, $V_{4}=\{e_{04},v_{1}\}$,
$V_{5}=\{e_{12},e_{34}\} \cup \{v_0\}$ and $V_{6}=\{e_{23},e_{45}\}$, we have $\chi_d^{tt}({W_4})=6$.
%------------------------------------------------------------
\item $n=5$. By the contrary, let $\ell=6$. Then \texttt{Fact 5}  implies $2b_1+b_2\geq 5$, and by \texttt{Fact 1} we have $(V_1,\cdots,V_6)=(4,4,4,2,1,1)$. But by considering the proof of Lemma \ref{alpha (T(W_n))=lceil 2n/3 rceil} we know that all of the maximum independent sets in $T(W_5)$ are the five sets $\{e_{0i},v_{i+1},e_{(i+2)(i+3)},v_{i+5}\}$ for $1\leq i\leq 5$, which only two of them are disjoint. Thus $V_i \cap V_j\neq \emptyset$ for some $1\leq i< j\leq 3$, a contradiction.  So $\ell\geq 7$, and since $(V_1,\cdots,V_7)$ is a TDC of $T(W_{5})$ where $V_{1}=\{v_{1},v_3,e_{02},e_{45}\}$, $V_{2}=\{v_{2},e_{34},e_{05}\}$,  $V_{3}=\{e_{03},v_4\}$, $V_{4}=\{e_{01},e_{23}\}$, $V_{5}=\{e_{04},v_5\}$, $V_{6}=\{v_0,e_{15}\}$, $V_{7}=\{e_{12}\}$, we have $\chi_d^{tt}({W_5})=7$.
%------------------------------------------------------------
\item $n=6$. By the contrary, let $\ell=7$. Then \texttt{Fact 5} implies $6\leq 2b_1+b_2\leq 7$. Since obviousely $b_{1}\geq 4$ implies $|V_1| > \alpha=4$, we assume  $b_{1} \leq 3$, and so $(b_1,b_2)=(3,0)$, $(3,1)$, $(3,2)$, $(3,3)$, $(3,4)$, $(2,2)$, $(2,3)$, $(2,4)$, $(2,5)$, $(1,4)$, $(1,5)$, $(1,6)$,  $(0,7)$. Since $(b_1,b_2)=(3,4)$, $(2,5)$, $(1,6)$, $(0,7)$ imply $\sum_{i=1}^{7}|V_i|\neq 3n+1$, and $(b_1,b_2)=(3,1)$, $(3,2)$, $(3,3)$, $(2,2)$, $(2,3)$, $(2,4)$, $(1,4)$, $(1,5)$ imply $|V_1|>\alpha=4$, which contradict \texttt{Fact 1}, we assume $(b_1,b_2)=(3,0)$. But this implies $(|V_1|,\cdots,|V_{7}|)=(4,4,4,4,1,1,1)$, by \texttt{Fact 1}, which is not possible. Because, by considering the proof of Lemma \ref{alpha (T(W_n))=lceil 2n/3 rceil}, the number of disjoint maximum independent sets in $T(W_6)$ is at most three. Therefore $\ell \geq 8$, and since  the coloring function $(V_1,\cdots,V_8)$ is a TDC of $T(W_6)$ where
$V_1=\{e_{12}, e_{34}, e_{56}, v_{0}\}$, $V_2=\{e_{23}, e_{45}, e_{16}\}$, $V_3=\{e_{01},v_{2}\}$, $V_4=\{e_{02},v_{3}\}$, $V_5=\{e_{03},v_{4}\}$, $V_6=\{e_{04},v_{5}\}$, $V_7=\{e_{05},v_{6}\}$, $V_8=\{e_{06},v_{1}\}$, we have $\chi_d^{tt}({W_6})=8$. 

%-------------------------------------------
\item $n=7$. By the contrary, let $\ell=8$. Then \texttt{Fact 5} implies $7\leq 2b_1+b_2\leq 8$. Since obviousely $b_{1}\geq 5$ implies $|V_1| > \alpha=5$, we assume  $b_{1} \leq 4$, and so $(b_1,b_2)=(4,0)$, $(4,1)$, $(4,2)$, $(4,3)$, $(4,4)$, $(3,1)$, $(3,2)$, $(3,3)$, $(3,4)$, $(3,5)$, $(2,3)$, $(2,4)$, $(2,5)$, $(2,6)$, $(1,5)$, $(1,6)$, $(1,7)$,  $(0,8)$. Since $(b_1,b_2)=(4,4)$, $(3,5)$, $(2,6)$, $(1,7)$, $(0,8)$ imply $\sum_{i=1}^{8}|V_i|\neq 3n+1$, and $(b_1,b_2)=(4,1)$, $(4,2)$, $(4,3)$, $(3,3)$, $(3,4)$, $(2,4)$, $(2,5)$, $(1,5)$, $(1,6)$ imply $|V_1|>\alpha=5$, which contradict \texttt{Fact 1}, we assume $(b_1,b_2)=(4,0)$, $(3,1)$, $(3,2)$, $(2,3)$. But then we have $4 \leq |V_3|\leq |V_{2}|\leq |V_1|\leq 5$, which is not possible. Because, by considering the proof of Lemma \ref{alpha (T(W_n))=lceil 2n/3 rceil}, the number of disjoint  independent sets of cardinalities four or five in $T(W_6)$ is at most two. So $\ell \geq 9$, and since the coloring function $(V_1,\cdots,V_9)$ is a TDC of $T(W_{7})$ where $V_{1}=\{e_{01},e_{34},e_{56},v_{2},v_{7}\}$, $V_{3}=\{e_{12},e_{45},e_{67},e_{03}\}$, $V_{5}=\{e_{23},e_{05},v_{1},v_{4},v_{6}\}$, $V_{7}=\{e_{17},v_{3},v_{5}\}$, $V_{2i}=\{e_{0(2i)}\}\mbox{ (for }1\leq i \leq 3)$, $V_{8}=\{v_{0}\}$, $V_{9}=\{e_{07}\}$, we have $\chi_d^{tt}({W_7})= 9$.
\end{itemize}
%----------------------------------
Figure \ref{TDTC8} shows $(\{v_{1}, v_{3}, e_{02}, e_{45}\},\{v_{2}, e_{34}, e_{05}\},\{v_{0}, e_{15}\},\{v_{4},e_{03}\}, \{v_{5},e_{04}\}, \{e_{01},e_{23}\}, \{e_{12}\})$ as a min-TDTC of $W_{5}$ (left) and as a min-TDC of $T(W_{5})$ (right).
%---------------------------
\begin{figure}[ht]
\centerline{\includegraphics[width=11cm, height=5cm]{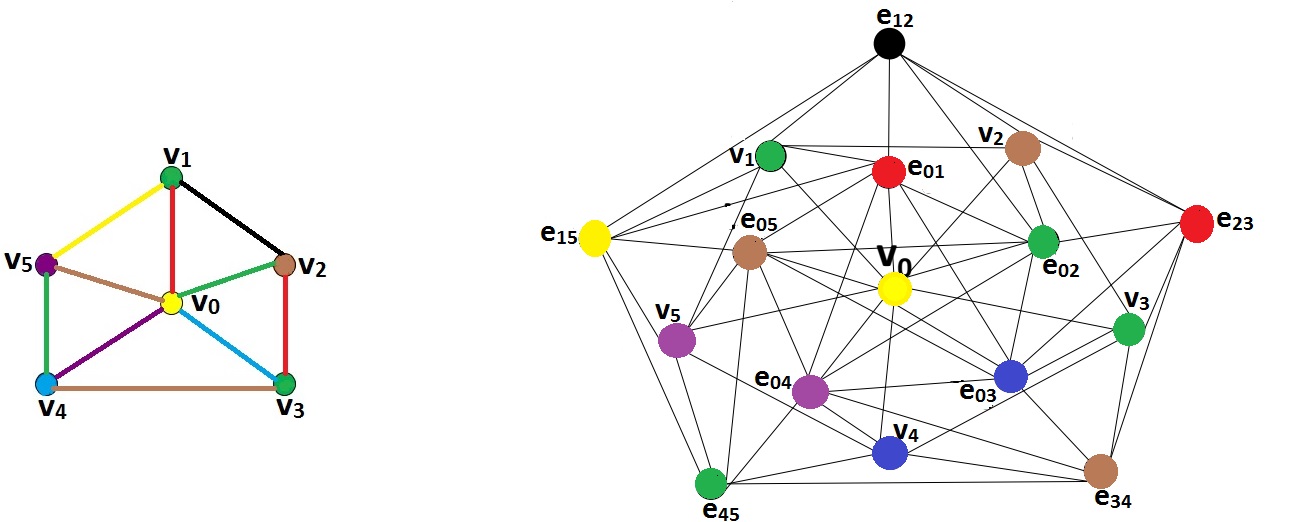}}
\vspace*{-0.25cm}
\caption{ A min-TDTC of $W_{5}$ (left)  and its corresponding  min-TDC of $T(W_{5})$ (right).}\label{TDTC8}
\end{figure}
%-----------------------------------------------------------------------

\end{proof}
%---------------------------------------
Proposition \ref{chi_d^t(T({W_n}))} shows that the upper bound given in Proposition \ref{chi_d^t =<g_t+min{chi(G-S)}} is tight for wheels of order more than 8.
%------------------------------------------

\section{Complete bipartite graphs}

Here, we calculate the total dominator total chromatic number of a complete bipartite graph $K_{m,n}=(V\cup U,E)$ in which $V\cup U$ is the partition of its vertex set to the independent sets $V=\{v_{i}: 1\leq i \leq m\}$, $U=\{u_{j}: 1\leq j \leq n\}$ and $E=\{v_iu_j~|~1\leq i \leq m, ~1\leq j \leq n \}$ is its edge set. 

%----------------------------------------------------------------------
\begin{prop}
\label{chi^{t}_d (T(K_{m,n}))}
For any complete bipartite graph $K_{m,n}$ in which $n\geq m\geq 1$, 
\begin{equation*}
\chi^{tt}_d (K_{m,n})=\left\{
\begin{array}{ll}
m+n  & \mbox{if }m=1,2~ \mbox{and}~ (m,n)\neq(1,1) \\
m+n+1  & \mbox{if }m\geq3~ \mbox{or}~ (m,n)=(1,1).
\end{array}
\right.
\end{equation*}
\end{prop}
\begin{proof}
Let $K_{m,n}=(V\cup U,E)$ be the descriptive complete bipartite graph in above of order $n+m\geq 2$. Hence $V\cup U\cup \mathcal{E}$ is a partition of the vertex set $T(K_{m,n})$ where $\mathcal{E}=\{e_{ij}~|~ 1\leq i \leq m, ~1\leq j \leq n  \}$. Since $T(K_{1,1})\cong K_{3}$ implies $\chi^{tt}_{d}(K_{1,1})=3$, we assume $n>m=1$. Let $f=(V_{1},V_{2},\cdots,V_{\ell})$ be a min-TDC of $T(K_{m,n})$. Since the subgraph of $T(K_{m,n})$ induced by $\{v_{1},e_{11},\cdots, e_{1n}\}$ is a complete graph of order $n+1$, we have $\chi^{t}_{d}(T(K_{m,n})\geq n+1$. As we have shown for $T(K_{1,3})$ in Figure \ref{TDTC3}, since $(V_{1},V_{2},\cdots,V_{n+1})$ is a TDC of $T(K_{1,n})$ where $V_{1}=\{e_{11}, u_{n}\}$, $V_{i}=\{e_{1i}, u_{i-1}\}$ for $2\leq i \leq n$, $V_{n+1}=\{ v_{1}\}$, which implies $\chi^{tt}_{d}(K_{1,n})= n+1$, we continue our proof in the following two cases.

\vspace{0.3cm}

\textbf{Case 1.} $n \geq m=2$. Let $\chi^{tt}_{d}(K_{m,n})= n+1$, and let $\mathcal{E}_{i}=\{ e_{ij}~|~ 1 \leq j \leq n \}$ for $i=1,2$. Since $T(K_{m,n})[\mathcal{E}_{1}] \cong T(K_{m,n})[\mathcal{E}_{2}] \cong K_{n}$ we have to color the vertices in $\mathcal{E}_{1}$ (and also in $\mathcal{E}_{2}$) by $n$ different colors. On the other hand, since $T(K_{m,n})[\mathcal{E}_{1} \cup  \mathcal{E}_{2}] \cong K_{n} \square K_{2}$ we conclude that $e_{1j}$ and $e_{2j}$ are not in a same color class for $1 \leq j \leq n$. Without loss of generality, we may assume $e_{1j} \in V_{j}$ for $1 \leq j \leq n$ and $ v_{1} \in V_{n+1}$. If $f(\mathcal{E}_{2})=\{1,2, \cdots,n\}$, then $v_{1} \nsucc_{t} V_{k}$ for each $1 \leq k \leq n$, because $N_{T(K_{m,n})}(v_{1}) \cap \mathcal{E}_{2}=\emptyset $ and $|V_{k}| \geq 2$ for each $1 \leq k \leq n$. So $n+1 \in f(\mathcal{E}_{2})$, and a color, say $1$, is not in $f(\mathcal{E}_{2})$. This implies $f(v_{2})=1$ and so $v_{1} \nsucc_{t} V_{k}$ for each $1 \leq k \leq n+1$. Thus $\ell \geq n+2=n+m$. Now by assumptions $V_{1}=\{e_{11},e_{2n}\}$, $V_{i}=\{e_{1i},e_{2(i-1)}\}$ for $2\leq i \leq n$, $V_{n+1}=V$, $V_{n+2}=U$, since the coloring function $(V_{1},V_{2},\cdots,V_{n+2})$ is a TDC of $T(K_{m,n})$, we have $\chi^{tt}_{d}(K_{m,n})= n+2$.

\vspace{0.3cm}

\textbf{Case 2.} $n \geq m \geq 3$. For $1\leq i \leq m$ let $\mathcal{E}_{i}=\{ e_{ij}~|~ 1 \leq j \leq n \}$, and for $ 1\leq j \leq n$ let $\mathcal{E}'_j=\{ e_{ij}~|~ 1 \leq i \leq m \}$. It can be easily seen that $T(K_{m,n})[\mathcal{E}_{i}] \cong K_{n}$, $T(K_{m,n})[\mathcal{E}'_{j}] \cong K_{m}$, $T(K_{m,n})[\mathcal{E}_{1} \cup \mathcal{E}_{2}\cup \cdots \cup \mathcal{E}_{m}] \cong T(K_{m,n})[\mathcal{E}'_{1} \cup \mathcal{E}'_{2}\cup \cdots \cup \mathcal{E}'_{n}]\cong K_{n} \square K_{m}$ and $T(K_{m,n})[\mathcal{E}_{i} \cup \{v_{i}\}] \cong K_{n+1}$, $T(K_{m,n})[\mathcal{E}'_{j} \cup \{u_{j}\}] \cong K_{m+1}$. By proving $\ell \geq n+m+1$ in the following two subcases, and by considering this fact that the coloring function $g$ with the criterion 
\begin{equation*}
\begin{array}{ll}
g(e_{ij}) \equiv j-i+1 \pmod{n} & \mbox{if } 1\leq i \leq m \mbox{ and } 1\leq j \leq n,\\
~~~~~~~~~~~g(v_{i})=n+i & \mbox{if } 1\leq i \leq m, \mbox{ and }\\
~~~~~~~~~~~ g(u_{i})=n+m+1 & \mbox{if } 1\leq i \leq n,
\end{array}
\end{equation*}
%----------------------
is a TDC of $T(K_{m,n})$ with $m+n+1$ color classes, we have $\chi^{tt}_{d}(K_{m,n})=m+n+1$.
%-----------------------------------

\begin{itemize}
\item{} $f(\mathcal{E}_{1} \cup \mathcal{E}_{2}\cup \cdots \cup \mathcal{E}_{m})=\{1,2, \ldots, n\}$. %( Notice this can be happend, because the function $f$ with the criterion $f(t_{ij}) \equiv j-i+1 \pmod{n} $ when $1 \leq i \leq m$ and $1 \leq j \leq n$ is a n-coloring of it.)
Since for each $1 \leq i \leq m$,  $v_{i} \succ_{t} V_{k_i}$ implies $f(V_{k_i}) \cap \{1,2, \ldots, n\}=\emptyset$ (because every color $1\leq i \leq n$ appears $m\geq 2$ times) and $V_{k_i} \subseteq U$, we have $\ell \geq n+1$. On the other hand, we see that for each $1\leq i \leq m$ and $1\leq j \leq n$, $e_{ij} \succ_{t} V_{k_{ij}}$ implies $V_{k_{ij}}\subset \{v_{i},u_{j}\}$. By the minimality of $f$, $n \geq m$ implies $V_{k_{ij}}=\{v_{i}\}$ for each $1 \leq i \leq m$. Now  the fact $f(V) \cap f(U)=\emptyset$ implies $\ell \geq n+m+1$. 

\item{} $\{1,2, \cdots, n+1\} \subseteq f(\mathcal{E}_{1} \cup \mathcal{E}_{2}\cup \cdots \cup \mathcal{E}_{m})$. We assume the min-TDC $f$ of $T(K_{m,n})$ is \emph{best} in this meaning that for every min-TDC $g$ of $T(K_{m,n})$, $|f(\mathcal{E}_{1} \cup \mathcal{E}_{2}\cup \cdots \cup \mathcal{E}_{m})| \leq |g(\mathcal{E}_{1} \cup \mathcal{E}_{2}\cup \cdots \cup \mathcal{E}_{m})|$. Then for each $1 \leq i \leq m$, $v_{i} \succ_{t} V_{k_i}$ implies $V_{k_i} \subseteq U \cup \mathcal{E}_{i}$ and specially if also $e_{ij} \in V_{k_i}$ for some $1 \leq i \leq n$, then $u_{j} \notin V_{k_i}$ and $f(e_{ij}) \notin f(\mathcal{E}_{1} \cup \mathcal{E}_{2}\cup \cdots \cup \mathcal{E}_{m})-f(\mathcal{E}_{i})$, that is, the color of $e_{ij}$ does not appear in the other vertices of $\mathcal{E}_{1} \cup \mathcal{E}_{2}\cup \cdots \cup \mathcal{E}_{m}-\mathcal{E}_i$. If every color in $f(\mathcal{E}_{1} \cup \mathcal{E}_{2}\cup \cdots \cup \mathcal{E}_{m})$ is appear at least two times, then similar to Case 1, we can prove that at least $m+1$ new color are needed for coloring of $V \cup U$, which implies $\ell \geq n+m+1$. Therefore, we assume there exists at least one color which is used for coloring of only one vertex in $\mathcal{E}_{1} \cup \mathcal{E}_{2}\cup \cdots \cup \mathcal{E}_{m}$. For $1 \leq i \leq m$ let $r_i$ be the number of colors which are used only for coloring of one vertex from $\mathcal{E}_{i}-(\mathcal{E}_{1}\cup \cdots \cup \mathcal{E}_{i-1})$. Without loss of generality, we may assume $r_1\geq r_2\geq \cdots \geq r_m$. We know $|f(\mathcal{E}_i)|= n$ for each $i$. Since $|f(\mathcal{E}_2)\cap f(\mathcal{E}_1)|\leq n-r_1$, we have $|f(\mathcal{E}_2)- f(\mathcal{E}_1)|\geq r_1$. In a similar way, we have $|f(\mathcal{E}_k)- \cup_{i=1}^{k-1}f(\mathcal{E}_i)|\geq \sum_{i=1}^{k-1} r_i$ for $3\leq k\leq m$. By summing this inequalities, we obtain 
%--------------------------------------------
\begin{equation}
\label{|f(T_1 ... T_m)|>= n+(m-1)r_1+(m-2)r_2+...+r_{m-1}}
|f(\mathcal{E}_{1} \cup \mathcal{E}_{2}\cup \cdots \cup \mathcal{E}_{m})|\geq n+(m-1)r_1+(m-2)r_2+\cdots +r_{m-1}.
\end{equation}
%-------------------------------------
Since (\ref{|f(T_1 ... T_m)|>= n+(m-1)r_1+(m-2)r_2+...+r_{m-1}}) implies $\ell\geq n+m+1$ when $r_1\geq 2$, we assume $r_1=1$. If $r_1=r_2=r_3=1$, then $m\geq 4$ and again (\ref{|f(T_1 ... T_m)|>= n+(m-1)r_1+(m-2)r_2+...+r_{m-1}}) implies $\ell\geq n+m+1$. Otherwise, there exists at least a vertex $e_{ij} \in \mathcal{E}_{i}$ for some $3 \leq i \leq m$ such that if $e_{ij} \succ_{t} V_{k_{ij}}$, then $V_{k_{ij}} \cap f(\mathcal{E}_{1} \cup \mathcal{E}_{2}\cup \cdots \cup \mathcal{E}_{m})=\emptyset$, that is, at least a new color  is needed, and  $V_{k_{ij}}\subset \{v_{i},u_{j}\}$. Since $f(V) \cap f(U)=\emptyset$, $V_{k_{ij}}=\{v_{i}\}$ implies that at least one new color is used to color some vertex in $U$, and similarly $V_{k_{ij}}=\{u_{j}\}$ implies that at least one new color is used to color some vertex in $V$. Therefore (\ref{|f(T_1 ... T_m)|>= n+(m-1)r_1+(m-2)r_2+...+r_{m-1}}) implies $\ell\geq (n+m-1)+1+1\geq n+m+1$. 
%there exists at least one new color so $\ell \geq n+m+1$. Otherwise, $|V_{n+1}| \geq 2$, $t_{11} \in V_{n+1}$ and $n+1 \in f(V)$ (or $t_{11} \in V_{n+1}$ and $n+1 \in f(U)$). Similary, we can proved, $\ell \geq n+m+1$. Since coloring function $g$ defined in case1, is a TDC of $T(K_{m,n})$ with $m+n+1$ color classes, we have $\chi^{t}_{d}(T(K_{m,n})= m+n+1$.
\end{itemize}
%------------------------------------
\begin{figure}[ht]
\centerline{\includegraphics[width=8cm, height=3cm]{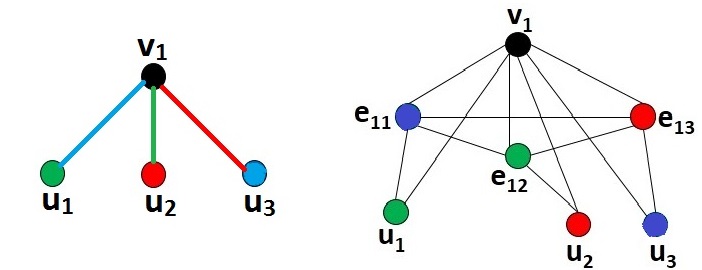}}
\vspace*{-0.25cm}
\caption{ A min-TDTC of $K_{1,3}$ (left) and its corresponding min-TDC of $T(K_{1,3})$ (right).}\label{TDTC33}
\end{figure}
%------------------------------------------

\end{proof}
%--------------------------------------------------------
As a result of Propositions \ref{chi_d^t(T({W_n}))} and \ref{chi^{t}_d (T(K_{m,n}))}, we have the following result.
%----------------------------------------------
\begin{thm}
For any $n\geq 3$, there exists a graph $G$ of order $n$ with $\chi^{tt}_d (G)=n$.
\end{thm}

\begin{proof}
Let $G=W_n$ when $n \geq 8$ or let  $G=K_{1,q}$ or $K_{2,q}$ of order at least 3. Then $\chi^{tt}_d (G)=n$ by Propositions \ref{chi_d^t(T({W_n}))} and \ref{chi^{t}_d (T(K_{m,n}))}.
\end{proof}
%----------------------------------------------------------------------------------------

\section{\bf Complete graphs}

\begin{figure}[ht]
\centerline{\includegraphics[width=7.5cm, height=5cm]{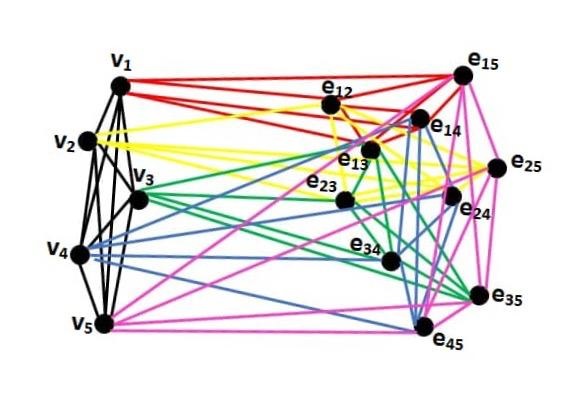}}
\vspace*{-0.6cm}
\caption{$T(K_{5})$ and its six edge-disjoint copies of $K_{5}$.}\label{T(K5)}
\end{figure}

From \cite{KKM2019}, we know that 
%-----------
\begin{prop}\emph{\cite{KKM2019}}
\label{chi^{t}_{d} (T(K_n)), n geq 8}
For any complete graph $K_n$ of order $n\geq 2$, $\chi^{tt}_{d}(K_n) \leq \lceil \frac{5n}{3} \rceil$.
\end{prop}
%------------------------------------------------------------------
Here, we show that the upper bound in Propsition \ref{chi^{t}_{d} (T(K_n)), n geq 8} is tight when $11\neq n\geq 9$. First we clarify more details on the total of a complete graph in the next observation. To more underestanding the observation, we have shown $T(K_5)$ in Figure \ref{T(K5)} as an example.
%---------------------------
%-------------------------------------------
\begin{obs}
\label{Obs on T(K_n)}
Let $T(K_n)$ be the total of a complete graph $K_n$ of order $n\geq 2$ with the vertex set $V=\{v_i~|~1\leq i \leq n\}$. Then $T(K_n)$ is $2(n-1)$-regular with the following properties.
\begin{itemize}
%-------------------
\item[1.] $T(K_n)=K_n^{v_0}\cup K_n^{v_1}\cup \cdots \cup K_n^{v_n}$ is  the partition of $T(K_n)$ to $n+1$ edge-disjoint copies of $K_n$ where $K_n^{v_0}=K_n$ and $V(K_n^{v_i})=\{v_i\}\cup \{e_{ij}~|~1\leq j\neq i \leq n\} $  for $1\leq i \leq n$.
%-------------------
\item[2.] $L(K_n)=T(K_n)-K_n=(K_n^{v_1}-\{v_1\})\cup \cdots \cup (K_n^{v_n}-\{v_n\})$ is  the partition of the line graph of $K_n$ to $n$ edge-disjoint copies of $K_{n-1}$.
%-------------------
\item[3.] $V(K_n^{v_i}) \cap V(K_n^{v_j})=\{e_{ij}\}$ for each $1\leq i< j \leq n$.
%-------------------
\item[4.] $V(K_n^{v_i}) \cap V(K_n^{v_0})=\{v_{i}\}$ for each $1\leq i \leq n$.
%-------------------
\item[5.] For every $x\in V(T(K_n))$, $N(x)=V(K_n^{v_i}) \cup V(K_n^{v_j})-\{x\}$ for some $0\leq i<j \leq n$.
%-----------------------------
\item[6.] For each $1\leq i \leq n$, the function $\phi_{i}$ on $V(T(K_n))$ with the criterion
\begin{equation*}
\phi_{i}(x)=\left\{
\begin{array}{ll}
v_i  & \mbox {if } x=v_i,\\
v_j  & \mbox {if } x=e_{ij},\\
e_{ij}  & \mbox {if } x=v_j,\\
x  & \mbox {otherwise}\\
\end{array}
\right.
\end{equation*}
is an authomorsim of $T(K_n)$ which replaces $K_n$ with $K_n^{v_i}$.  And so $\phi_j \circ \phi_i^{-1}$ is an authomorsim of $T(K_n)$ which replaces $K_n^{v_i}$ with $K_n^{v_j}$.
%-------------------
%\item[7.] For each $1\leq i\neq j \leq n$, the function $\varphi_{ij}$ on $V(T(K_n))$ with the criterion
%\begin{equation*}\varphi_{ij}(x)=\left\{
%\begin{array}{ll}v_j  & \mbox {if } x=v_i,\\v_i  & \mbox {if } x=v_j,\\e_{jk}  & \mbox {if } x=e_{ik} \mbox{ for } 1\leq i\neq k\neq j \leq n,\\e_{ik}  & \mbox {if } x=e_{jk} \mbox{ for } 1\leq i\neq k\neq j \leq n,\\x  & \mbox {otherwise}\\\end{array}\right.\end{equation*}is an authomorsim of $T(K_n)$ which replace $K_n^{v_i}$ with $K_n^{v_j}$.
%-------------------
\end{itemize} 
\end{obs}
%-------------------------------------------------------
The following proposition from \cite{KK2017} is useful for our investigation.
%--------------------------------------------
%\begin{prop} \emph{\cite{KK2017}}\label{alpha_{mix}(K_n)=lceil n/2 rceil}any complete graph $K_n$ of order $n\geq 2$, $\alpha_{mix}(K_n)=\lceil n/2 \rceil$.\end{prop}
%------------------------------------------------------
\begin{prop} \emph{\cite{KK2017}}
\label{gamma_t T(K_n)}
For any complete graph $K_n$ of order $n\geq 2$, 
\begin{itemize}
\item[1.] $\gamma_{tm}(K_n)=\gamma_t(T(K_n))=\lceil \frac{5n}{3}\rceil  -n$,

\item[2.] $\alpha_{mix}(K_n)=\alpha(T(K_n))=\lceil \frac{n}{2}  \rceil$.
\end{itemize}
\end{prop}
%----------------------------------------------------------
\begin{prop} 
\label{chi^{t}_{d} (T(K_n)) 2=<n=<7}
For any complete graph $K_n$ of order $n \geq 2$,
\begin{equation*}
\chi^{tt}_{d}(K_n)=\left\{
\begin{array}{ll}
\lceil \frac{5n}{3}  \rceil -2  & \mbox {if } n=3,4,5,\\
\lceil \frac{5n}{3}  \rceil -1  & \mbox {if } n=2,6,7,8,11\\
\lceil \frac{5n}{3}  \rceil   & \mbox {if } n\geq 9 ~\mbox {and } n\neq 11.
\end{array}
\right.
\end{equation*}
\end{prop}

\begin{proof}
Let $ K_{n} $ be a complete graph with the vertex set $V=\{v_i~|~1\leq i \leq n\}$. By Propositions \ref{chi^{t}_d T(P_n)}, \ref{chi^{t}_d T(C_n)}, \ref{chi_d^t(T({W_n}))}, we may assume $n\geq 5$.  Then $V(T(K_{n}))=V\cup \mathcal{E}$ where $\mathcal{E}=\{e_{i  j}~|~ 1\leq i<j \leq n \}$. Let  $f=(V_{1},V_{2},\cdots,V_{\ell})$ be a min-TDC of $T(K_{n})$ in which $|V_1|\geq |V_2|\geq \cdots \geq |V_{\ell}|$, and for $1 \leq i \leq \alpha:=\lceil \frac{n}{2} \rceil $ (recall from Proposition \ref{gamma_t T(K_n)}(1) that $\lceil \frac{n}{2} \rceil =\alpha_{mix}(K_n)$) let $\mathcal{A}_{i}=\{V_{k}~|~ v \succ_{t} V_{k}\mbox{ and }|V_{k}| = i \mbox{ for some }v \in V\cup \mathcal{E}\}$ and $|\mathcal{A}_{i}|=a_i$. By the next nine facts, we continue our proof in the following two cases.
%------------------------
\begin{itemize}
\item[$\star$] \texttt{Fact 1.} %$V_1\cup \cdots \cup V_{\ell}=V\cup \mathcal{E}=CN(V_1)\cup \cdots \cup CN(V_{\ell})$, which imply 
$\sum_{i=1}^{\ell}|V_i|=\frac{n(n+1)}{2}$ and $\sum_{i=1}^{\ell}|CN(V_i)|\geq \frac{n(n+1)}{2}$ by (\ref{|V|=sum_{i=1}^{ell} |V_i|}) and (\ref{CN(V_i)=V}), respectively. Also $|V_{k}|\leq \lceil \frac{n}{2}\rceil $ for each $1\leq k\leq \ell$.
 
\item[$\star$] \texttt{Fact 2.} For any $v \in V\cup \mathcal{E}$, if $v \succ_{t} V_{k}$ for some $1\leq k\leq \ell$, then $|V_{k}|\leq 2$. Because $N(v)=V(K_n^{v_i}) \cup V(K_n^{v_j})-\{v\}$, for some $0\leq i<j \leq n$ (by Observation \ref{Obs on T(K_n)}(5)), implies $|V_k\cap V(K_n^{v_i})|\leq 1 $ and $|V_k\cap V(K_n^{v_j})|\leq 1 $.

\item[$\star$] \texttt{Fact 3.} If $V_k=\{v_i,e_{pq}\}$ for different indices $i$, $p$, $q$, then $CN(V_k)=\{v_p, v_q, e_{ip}, e_{iq}\}$, and if $V_k=\{e_{rs},e_{pq}\}$ for different indices $p$, $q$, $r$, $s$, then $CN(V_k)=\{e_{rp},e_{rq},e_{sp}, e_{sq}\}$.

\item [$\star$] \texttt{Fact 4.} If $V_k=\{v_i\}$ for some $i$, then $CN(V_k)=V \cup \{e_{ij}~|~1\leq j\neq i \leq n\}-\{v_i\}$, and if $V_k=\{e_{pq}\}$ for some $p\neq q$, then $CN(V_k)=\{e_{ij}~|~|\{p,q\}\cap \{i,j\}|=1\}\cup \{v_p,v_q\}$.

\item[$\star$] \texttt{Fact 5.} $(2n-2)a_1+4a_2\geq \frac{n(n+1)}{2}$.\\ Because
\begin{equation*}
\begin{array}{llll}
\frac{n(n+1)}{2} &   =   & |V\cup \mathcal{E}| &\\
                          & \leq  & \Sigma_{|V_k|\leq 2}|CN(V_k)| & \mbox{(by \texttt{Fact 2})}\\
                          &  =    & \Sigma_{|V_k|=1}|CN(V_k)|+ \Sigma_{|V_k|=2}|CN(V_k)|& \\
                          & \leq & (2n-2)a_1+4a_2. &
     \end{array}
%\right.
\end{equation*} 
%--------------
\item[$\star$] \texttt{Fact 6.} $\lceil \frac{5n}{3}\rceil  -n \leq a_1+a_2\leq \ell $.\\ Because the set $S$ with this property that $|S\cap V_i|=1$ for each $V_i\in \mathcal{A}_1\cup \mathcal{A}_2$ is a TDS of $T(K_n)$ (by \texttt{Fact 2} and Proposition \ref{gamma_t T(K_n)}(2) for left), and $a_1+\cdots +a_{\alpha}=\ell$ (for right).

\item[$\star$] \texttt{Fact 7.} $\lceil \frac{n(n+1)/2 -4(\lceil 5n/3\rceil  -n)}{2n-6} \rceil \leq a_1 \leq \lfloor \frac{\alpha \ell -n(n+1)/2}{\alpha -1} \rfloor $.\\
Because  the lower bound can be obtained by \texttt{Facts 5,6}, and the upper bound can be obtained by
%-------------------------------------
\begin{equation*}
%\label{2n-1-a_1 leq (ell-a_1)alpha}
\begin{array}{llll}
\frac{n(n+1)}{2}-a_1             &   =   & |V(T(K_n))|-|\mathcal{A}_1| \\
                          & =     & \Sigma_{|V_i|\geq 2} |V_i|\\
                          & \leq & (\ell-a_1)\alpha.
     \end{array}
\end{equation*} 
%------------------------------------------
\item[$\star$] \texttt{Fact 8.} $\lceil \frac{(2n-2)(\lceil 5n/3\rceil  -n)-n(n+1)/2}{2n-6} \rceil \leq a_2 \leq \ell-a_{1} $ (by \texttt{Facts 5,6,7}).

\item[$\star$] \texttt{Fact 9.} For any $V_i=\{e_{rs}\},V_j=\{e_{pq}\}\in \mathcal{A}_1$ (it is allowed $r=s$ or $p=q$, and in this case $e_{ii}$ is the same $v_i$),
\begin{equation*}
|CN(V_i)\cap CN(V_j)|=\left\{
\begin{array}{ll}
4     & \mbox {if } \{r,s\}\cap \{p,q\}=\emptyset, \\
n-1  & \mbox {if } \{r,s\}\cap \{p,q\}\neq \emptyset.
\end{array}
\right.
\end{equation*}

\end{itemize}

%--------------------------------------------------------------------
\textbf{Case 1.} $5 \leq n \leq 8$ or $n=11$.
\begin{itemize}
\item $n=5$. Let $\ell=6$. Then $(a_1,a_2)=(0,5)$, $(0,6)$, $(1,5)$. Because $0 \leq  a_1 \leq 1$ and $5\leq a_2 \leq 6$ by \texttt{Facts 7, 8}. Since $(a_1,a_2)= (0,6)$, $(1,5)$ imply $\sum_{i=1}^{6}|V_i|\neq \frac{n(n+1)}{2}$, and $(a_1,a_2)=(0,5)$ implies $|V_1|>\alpha=3$, which contradict \texttt{Fact 1}, we have  $\ell\geq 7$. Now since $(V_{1}, \cdots, V_{7})$ is a TDC of $T(K_{5})$ where $V_{1}=\{v_{3},e_{12},e_{45}\}$, $V_{2}=\{v_{4},e_{23},e_{15}\}$, $V_{3}=\{v_{5}, e_{13},e_{24}\}$, $V_{4}=\{e_{25},e_{34}\}$, $V_{5}=\{e_{35},e_{14}\}$, $V_{6}=\{v_{1}\}$, $V_{7}=\{v_{2}\}$, we have $\chi^{tt}_{d}(K_5) =7=\lceil \frac{5n}{3}  \rceil -2$.
%----------------------------------------------------------------------
\item $n=6$. Let $\ell=8$. Then $(a_1,a_2)=(1,4)$, $(1,5)$, $(1,6)$, $(1,7)$. Because $ a_1 =1$ and $4\leq a_2 \leq 7$ by \texttt{Facts 7, 8}. Since $(a_1,a_2)= (1,7)$ implies $\sum_{i=1}^{8}|V_i|\neq \frac{n(n+1)}{2}$, and $(a_1,a_2)=(1,4)$, $(1,5)$, $(1,6)$ imply $|V_1|>\alpha=3$, which contradict \texttt{Fact 1}, we have  $\ell\geq 9$. Now since $(V_{1}, \cdots, V_{9})$ is a TDC of $T(K_{6})$ where $V_{1}=\{v_{3},e_{12},e_{45}\}$, $V_{2}=\{v_{4},e_{13},e_{26}\}$, $V_{3}=\{v_{5},e_{16},e_{23}\}$, $V_{4}=\{v_{6},e_{14},e_{25}\}$, $V_{5}=\{e_{36},e_{15},e_{24}\}$, $V_{6}=\{ e_{34},e_{56}\}$,  $V_{7}=\{e_{35},e_{46}\}$, $V_{8}=\{v_{1}\}$, $V_{9}=\{v_{2}\}$, we have $\chi^{tt}_{d}(K_6) =9=\lceil \frac{5n}{3}  \rceil -1$. 
%-------------------------------------------------------------------
\item $n=7$. Let $\ell=10$. Then $(a_1,a_2)=(1,4)$, $(1,5)$, $(1,6)$, $(1,7)$, $(1,8)$, $(1,9)$, $(2,4)$, $(2,5)$, $(2,6)$, $(2,7)$, $(2,8)$, $(3,4)$, $(3,5)$, $(3,6)$, $(3,7)$, $(4,4)$, $(4,5)$, $(4,6)$. Because $1 \leq a_1 \leq 4$ and $4\leq a_2 \leq 10-a_1$ by \texttt{Facts 7, 8}.  Since $(a_1,a_2)=(1,9)$, $(2,8)$, $(3,7)$, $(4,6)$ imply $\sum_{i=1}^{10}|V_i|\neq \frac{n(n+1)}{2}$, and $(a_1,a_2)=(1,5)$, $(1,6)$, $(1,7)$, $(1,8)$, $(2,4)$, $(2,5)$, $(2,6)$, $(2,7)$, $(3,4)$, $(3,5)$, $(3,6)$, $(4,4)$, $(4,5)$ imply $|V_1|>\alpha=4$, which contradict \texttt{Fact 1}, we have $(a_1,a_2)=(1,4)$. %Then either $V_{10}=\{v_i\}$ for some $1\leq i \leq n$ or $V_{10}=\{e_{pq}\}$ for some $1\leq p<q \leq n$. 
By Observation \ref{Obs on T(K_n)}(6), we may assume $V_{10}=\{v_i\}$ for some $1\leq i \leq n$. Then $v_{i} \succ_{t} V_{k}$ for some $k\neq 10$ implies $V_{k}=\{v_{p}, e_{iq}\}$ for some three different indices $i$, $p$,  $q$. Since $|CN(V_k)\cup CN(V_{10})|=|CN(V_k)|+|CN(V_{10})|-|CN(V_k)\cap CN(V_{10})|=4+12-4=12$ (by \texttt{Facts 3, 4}), we reach to this contradiction that
\begin{equation*}
\begin{array}{llll}
28                      &  =  & \frac{n(n+1)}{2} &  \\
                          &  \leq  & \Sigma_{i=6}^{10}|CN(V_k)| &  \\
                          &  \leq  & \Sigma_{k\neq i=6}^{9}|CN(V_k)|+ \Sigma_{i=k,10}|CN(V_i)| & \\
                          & \leq & 3\times 4 +12 & (\mbox{because } |CN(V_i)|=4 \mbox{ by }\texttt{Fact 3}) \\
                          & = & 24.& 
     \end{array}
%\right.
\end{equation*} 
Thus $\ell\geq 11$, and since $(V_{1},\cdots,V_{11})$ is a TDC of $T(K_{7})$ where
%----------------
\begin{equation*}
\begin{array}{ll}
V_{1}=\{v_{4}, e_{16}, e_{25}, e_{37}\}, ~V_{2}=\{v_{5}, e_{67}, e_{13}, e_{24}\},~V_{3}=\{v_{6}, e_{15}, e_{23}, e_{47}\},\\
V_{4}=\{v_{7}, e_{35}, e_{26}, e_{14}\}, ~V_{5}=\{v_{3}, e_{46},e_{57}\}, ~~~~~~V_{6}=\{v_{1}, e_{27}, e_{36}\},\\
 V_{7}=\{v_{2}, e_{34}\},~V_{8}=\{ e_{12}\},~V_{9}=\{ e_{45}\},~V_{10}=\{ e_{56}\},~V_{11}=\{ e_{17}\},
\end{array}
\end{equation*}
%-------------------------
we have $\chi^{tt}_{d}(K_7) =11=\lceil \frac{5n}{3} \rceil -1$.
%------------------------------------------------------------------------------------------
\item $n=8$. Let $\ell=12$. Then $(a_1,a_2)=(2,5)$, $(2,6)$, $(2,7)$, $(2,8)$, $(2,9)$, $(2,10)$, $(3,5)$, $(3,6)$, $(3,7)$, $(3,8)$, $(3,9)$,  $(4,5)$, $(4,6)$, $(4,7)$, $(4,8)$.  Because $2 \leq a_1 \leq 4$ and $5\leq a_2 \leq 12-a_1$ by \texttt{Facts  7, 8}. Since $(a_1,a_2)=(2,10)$, $(3,9)$, $(4,8)$ imply $\sum_{i=1}^{12}|V_i|\neq \frac{n(n+1)}{2}$, and $(a_1,a_2)=(2,5)$,$(2,6)$, $(2,7)$, $(2,8)$, $(2,9)$,  $(3,5)$, $(3,6)$, $(3,7)$, $(3,8)$,  $(4,5)$, $(4,6)$, $(4,7)$  imply $|V_1|>\alpha=4$, which contradict \texttt{Fact 1},  we have $\ell\geq 13$. Now since $(V_{1},\cdots,V_{13})$ is a TDC of $T(K_{8})$ where 
\begin{equation*}
\begin{array}{ll}
 V_{1}=\{ v_{8} , e_{13}, e_{24},e_{56} \},~V_{2}=\{ v_{7} , e_{25}, e_{36},e_{48} \},~V_{3}=\{ v_{6} , e_{18}, e_{27},e_{34} \},\\
 V_{4}=\{ v_{5} , e_{16}, e_{28},e_{37} \},~ V_{5}=\{ v_{4} , e_{17}, e_{26},e_{35} \},~V_{6}=\{ v_{2} , e_{15}, e_{47},e_{38} \},\\
V_{7}=\{ e_{57}, e_{68}\},~V_{8}=\{ e_{46}, e_{58}\},~ V_{9}=\{ e_{45}, e_{67}\},~V_{10}=\{  e_{14}, e_{78}\},\\V_{11}=\{v_{3}, e_{12}\},~V_{12}=\{ e_{23}\},~V_{13}=\{v_{1}\},\\
\end{array}
\end{equation*}
 we have $\chi^{tt}_{d}(K_8) =13=\lceil \frac{5n}{3}  \rceil -1$. 
 %-------------------
 \item $n=11$. Let $\ell=\lceil \frac{5n}{3} \rceil-2=17$. Then $(a_1,a_2)=(3,6)$, $(3,7)$, $(3,8)$, $(3,9)$, $(3,10)$, $(3,11)$, $(3,12)$, $(3,13)$, $(3,14)$, $(4,6)$, $(4,7)$, $(4,8)$, $(4,9)$, $(4,10)$, $(4,11)$, $(4,12)$, $(4,13)$, $(5,6)$, $(5,7)$, $(5,8)$, $(5,9)$, $(5,10)$, $(5,11)$, $(5,12)$, $(6,6)$, $(6,7)$, $(6,8)$, $(6,9)$, $(6,10)$, $(6,11)$, $(7,6)$, $(7,7)$, $(7,8)$, $(7,9)$, $(7,10)$.  Because $3 \leq a_1 \leq 7$ and $6\leq a_2 \leq 17-a_1$ by \texttt{Facts  7, 8}. Since $\sum_{i=1}^{17}|V_i|\neq \frac{n(n+1)}{2}$ when $(a_1,a_2)=(3,14)$, $(4,13)$, $(5,12)$,$(6,11)$, $(7,10)$ and  $|V_1|>\alpha=6$ in the other cases, which contradict \texttt{Fact 1},  we have $\ell\geq 18$. Now since $(V_{1},\cdots,V_{18})$ is a TDC of $T(K_{11})$ where 
 %-------------------------------
\begin{equation*}
\begin{array}{ll}
V_{1}=\{ v_{11} , e_{1(10)}, e_{26},e_{37} ,e_{48} ,e_{59}\},~~~V_{2}=\{ v_{10} , e_{19}, e_{28},e_{35} ,e_{46} ,e_{7(11)}\},\\
V_{3}=\{ v_{9} , e_{1(11)}, e_{27},e_{34} ,e_{8(10)} ,e_{56}\},~V_{4}=\{ v_{8} , e_{14}, e_{2(11)},e_{39} ,e_{6(10)} ,e_{57}\},\\
V_{5}=\{ v_{7} , e_{18}, e_{29},e_{36} ,e_{4(10)} ,e_{5(11)}\},~V_{6}=\{ v_{4} , e_{15}, e_{2(10)},e_{3(11)} ,e_{68} ,e_{79}\},\\
V_{7}=\{ v_{5} , e_{16}, e_{24},e_{3(10)} ,e_{9(11)} \},~~~~~~V_{8}=\{ v_{6} , e_{17}, e_{25},e_{38} ,e_{49} \},\\
V_{9}=\{ v_{3} , e_{12}, e_{4(11)},e_{58} ,e_{7(10)} \},~~~~~~V_{10}=\{ e_{6(11)}, e_{9(10)}\},\\
V_{11}=\{ e_{47}, e_{5(10)}\},~V_{12}=\{ e_{89} , e_{(10)(11)} \},~V_{13}=\{ v_{2} , e_{13} \},~V_{14}=\{  e_{67}, e_{8(11)}\},\\
V_{15}=\{ e_{69}, e_{78}\},~V_{16}=\{v_{1}\},~V_{17}=\{ e_{23}\},~V_{18}=\{ e_{45}\},\\
\end{array}
\end{equation*}
 %-------------------------------
 we have $\chi^{tt}_{d}(K_{11}) =18=\lceil \frac{5n}{3}  \rceil -1$. 
 %-------------------
 \end{itemize}
 %-------------------------------
\textbf{Case 2. } $n\geq 9$ and $n\neq 11$. Let $\ell=\lceil \frac{5n}{3}  \rceil  -1$, and let
\begin{equation*}
\begin{array}{ll}
m_1=\lceil \frac{n(n+1)/2 -4(\lceil 5n/3\rceil  -n)}{2n-6} \rceil, ~~~~ M_1=\lfloor \frac{\alpha \ell -n(n+1)/2}{\alpha -1} \rfloor,\\
m_2=\lceil \frac{(2n-2)(\lceil 5n/3\rceil  -n)-n(n+1)/2}{2n-6} \rceil, ~~~~M_2=\ell-a_1.
\end{array}
\end{equation*}
Then $(a_1,a_2) =(x_1+i,x_2+j)$ for some $0\leq i \leq M_1-m_1$ and some $0\leq j \leq M_2-m_2$. In the following cases, we show that $\ell=\lceil \frac{5n}{3}  \rceil  -1$ leads us to a contrdiction.
%----------------------------------------------------------
\begin{itemize}

\item Either $n $ is even or $(a_1,a_2)\neq(m_{1},m_{2})$ when $n $ is odd.  Then, by \texttt{Fact 1}, we must  have $a_1+a_2 \leq \ell-1$ and $|V_1| \leq \lceil \frac{n}{2}\rceil$, and so by assumptions $a_1=m_1+i$ and $a_2=m_2+j$ for some $0 \leq i \leq M_{1}-m_{1}$ and some $0 \leq j \leq M_{2}-m_{2}$, we have
%--------------------
\begin{equation*}
\label{|V_1| leq rceil >lceil n/2rceil}
\begin{array}{llll}
 \sum_{i=1}^{\ell} |V_i| & = & \sum_{i=1}^{\ell-(a_1+a_2)} |V_i|+\sum_{i=\ell-(a_1+a_2)+1}^{\ell} |V_i|\\
& \leq & (\ell-a_1-a_2)\lceil \frac{n}{2} \rceil+a_1+2a_2\\
& =    & (\ell-m_1-m_2)\lceil \frac{n}{2} \rceil+(m_1+2m_2)+(1-\lceil \frac{n}{2} \rceil)i+(2-\lceil \frac{n}{2} \rceil)j\\ 
& \leq & (\ell-m_1-m_2)\lceil \frac{n}{2} \rceil+(m_1+2m_2)-(4i+3j)~~ ( \mbox{because } n\geq 9)\\ 
& \leq & (\ell-m_1-m_2-\epsilon )\lceil \frac{n}{2} \rceil+(m_1+2m_{2}+2\epsilon)\\
& <    &  \frac{n(n+1)}{2},
\end{array}
\end{equation*}
%-------------------------------
which contradicts \texttt{Fact 1} (where $\epsilon$ is 0 when $n$ is even and is 1 otherwise). 
%-----------------------------------------------

\item  $n \geq 13$ is odd and $(a_1,a_2)=(m_{1},m_{2})$. Then $a_1= \lfloor \frac{n+1}{4} \rfloor \geq 3 $ and 
%-------------------------
\begin{equation*}
a_2=\left\{
\begin{array}{ll}
\lceil \frac{5n}{12} \rceil   & \mbox {if } n \equiv 3 \pmod{12},\\
\lceil \frac{5n}{12} \rceil +1  & \mbox {if } n \not \equiv 3 \pmod{12}.
\end{array}
\right.
\end{equation*}
%-------------------------
Let $\{V_i~|~V_i\in \mathcal{A}_1\}=\{V_i~|~i\in I\}$ where $I=\{\ell-i ~|~|~0\leq i \leq \ell-a_1+1\}$. Let $z=\Sigma_{i,j\in I-\{t\}} |CN(V_i) \cap CN(V_i)|$ for some $t\in I$. Since $|CN(V_t) \cap CN(V_i)| \geq 4 $ for each $i\in I-\{t\}$ (by \texttt{Fact 9}) and $CN(V_t) \cap CN(V_i)\cap CN(V_j)=\emptyset$ for each 2-subset $\{i,j\} \subseteq I-\{t\}$, we conclude 
%----------------------------------
\begin{equation}
\label{Sigma_{i,j in I} |CN(V_i) cap CN(V_i)| geq & z+ {a_1-1 choose 1}}
\begin{array}{llll}
\Sigma_{i,j\in I} |CN(V_i) \cap CN(V_j)| & = & z+\Sigma_{i,j\in I-\{t\}} |CN(V_i) \cap CN(V_j)|\\
& \geq & z+4(a_1-1)\\
& =    & z+ 4{a_1-1 \choose 1}
\end{array}
\end{equation}
%----------------------------------
Since $z=4{3 \choose 2}$ when $a_1=3$, by induction on $a_1\geq 3$ and (\ref{Sigma_{i,j in I} |CN(V_i) cap CN(V_i)| geq & z+ {a_1-1 choose 1}}), we will have
%----------------------------------
\begin{equation*}
\label{Sigma_{i,j in I} |CN(V_i) cap CN(V_j)| geq & 4{a_1 choose 2}}
\begin{array}{llll}
\Sigma_{i,j\in I} |CN(V_i) \cap CN(V_j)| & \geq &  4{a_1-1 \choose 2}+4{a_1-1 \choose 1}\\
 & = &4 {a_1 \choose 2}.
\end{array}
\end{equation*}
%--------------------
So
\begin{equation*}
\label{|V_1| leq rceil >lceil n/2rceil}
\begin{array}{llll}
 \sum_{i=1}^{\ell} |CN(V_i)| & = & \sum_{V_i\in \mathcal{A}_2} |CN(V_i)|+\sum_{V_i\in \mathcal{A}_1} |CN(V_i)|\\
& \leq & 4a_2 +(2n-2)a_1-4 {a_1 \choose 2} \\
& <    & \frac{n(n+1)}{2},
\end{array}
\end{equation*}
%-------------------------------
which contradicts \texttt{Fact 1}.

\item  $n=9$ and $(a_1,a_2)=(2,5)$. Then $\ell=14$, $\mathcal{A}_1=\{V_{13}, V_{14}\}$ and $\mathcal{A}_2=\{V_{i}~|~8\leq i \leq 12\}$. If $T(K_{n}) [V_{13} \cup V_{14} ] \cong K_{2}$, then, by $|CN(V_{13}) \cup CN(V_{14})|=3n-3=24$ and \texttt{Facts 2, 3}, we have
%--------------------
\begin{equation*}
\label{|V_1| leq rceil >lceil n/2rceil}
\begin{array}{llll}
| \cup_{i=1}^{\ell} CN(V_i)| & = & | \cup_{i\in \mathcal{A}_1 \cup \mathcal{A}_2} CN(V_i)|\\
& \leq & | \cup_{i\in \mathcal{A}_1} CN(V_i)| + | \cup_{i\in \mathcal{A}_2} CN(V_i)|\\
& \leq & 24+4a_2\\
& = & 44\\
& <  & 45 \\ 
& =    &  \frac{n(n+1)}{2},
\end{array}
\end{equation*}
%-------------------------------
which contradicts \texttt{Fact 1}. So $V_{13} \cup V_{14}$ is an independent set, and by Obsevation \ref{Obs on T(K_n)}(6), we may assume $V_{13}=\{v_{1}\}$ and $V_{14}= \{e_{23}\}$, and so $|CN(V_{13}) \cup CN(V_{14})|=28$. Then the assumption $v_{1} \succ_{t} V_{12}$ implies $V_{12}=\{v_{p}, e_{1q}\}$ for some $p\neq q$ by Obsevation \ref{Obs on T(K_n)}(6). If $\{p,q\} \cap \{2,3\}=\emptyset$, then each of the six numbers 4, 5, 6, 7, 8, 9 must be appeared three times in the indices of the elements of $V_8\cup \cdots \cup V_{11}$, which is not possible. Because the number of indices in the elements of each $V_i\in \mathcal{A}_2$, is at most four. So, we have $\{p,q\} \cap \{2,3\}\neq \emptyset$. Then the number of appearing all of the six numbers 4, 5, 6, 7, 8, 9 (by allowing repeating numbers) as indices of the elements of $V_8\cup \cdots \cup V_{11}$ can be reduced to 16, but then $e_{23} \nsucc_{t} V_{k}$ for each $k$. 
%-------------------------------------------
\end{itemize}
%-------------------------------------------
Therefore $\ell\geq \lceil \frac{5n}{3} \rceil $, and in fact $\chi^{t}_d (T(K_{n}))=\lceil \frac{5n}{3}  \rceil $ by Proposition \ref{chi^{t}_{d} (T(K_n)), n geq 8}.
%----------------------------------------------
%DELET this figure. In Figure \ref{TDTC9}, $(\{v_{1}, e_{23}\},\{v_{2},e_{14}\}, \{v_{3},e_{24}\}, \{v_{4},e_{13}\} \{e_{12},e_{35}\})$  is a min-TDC of $T(K_{4})$.
%---------------------------
%\begin{figure}[ht]
%\centerline{\includegraphics[width=10cm, height=4cm]{TDTC99}}
%\vspace*{-0.25cm}
%\caption{ A min-TDTC of $K_{4}$ (left)  and its corresponding  min-TDC of $T(K_{4})$ (right).}\label{TDTC9}
%\end{figure}
\end{proof}
%----------------------------------------------------------------------
%---------REFERENCES -------------------------------------------

\end{document}